# Is the Halting probability a Dedekind real number?

Bhupinder Singh Anand

In a recent historical overview, Cristian S. Calude, Elena Calude, and Solomon Marcus identify eight stages in the development of the concept of a mathematical proof in support of an ambitious conjecture: we can express classical mathematical concepts adequately only in a mathematical language in which both truth and provability are essentially unverifiable. In this paper we show, firstly, that the concepts underlying their thesis can, however, be interpreted constructively; and, secondly, that an implicit thesis in the authors' arguments implies that the probability of a given Turing machine halting on a given input cannot be expressed as a Dedekind real number.

## Contents

## 1. Introduction

In a recently arXived paper [Ca01], "Passages of Proof", Cristian S. Calude, Elena Calude, and Solomon Marcus conjecture that:

> Knowledge is acquired through reason and by experiment. For a long time mathematical proofs required only reason; this might be no longer true.

As the basis for their belief, they identify eight stages in the development of the concept of a mathematical proof:

(*a*) The first period was that of pre-Greek mathematics, for instance the Babylonian one, dominated by observation, intuition and experience.

(*b*) The second period was started by Greeks such as Pythagoras and is characterised by the discovery of deductive mathematics, based on theorems.

(*c*) ... with Galilei, Descartes, Newton and Leibniz, the mathematical language became more and more a mixed language, characterized by a balance between its natural and artificial components. ... This was the third step in the development of mathematical proofs.

(*d*) The fourth step is associated with the so-called epsilon rigour, so important in mathematical analysis; it occurred in the XIXth century and it is associated with names such as A. Cauchy and K. Weierstrass.

(*e*) The fifth period begun with the end of the XIXth century, when Aristotle's logic, underlining mathematical proofs for two thousands years, entered a crisis with the challenge of the principle of non-contradiction.

(*f*) The sixth period begins with Godel's incompleteness theorem (1931), for many meaning the unavoidable failure of any attempt to formalise the whole mathematics.

(*g*) The seventh period belongs to the second half of the XXth century, when algorithmic proofs become acceptable only when their complexities were not too high.

(*h*) With the eighth stage, proofs are no longer exclusively based on logic and deduction, but also empirical and experimental factors.

## 1.1 What is proof?

The authors [Ca01] then ask, rhetorically:

> What is a mathematical proof ? At a first glance the answer seems obvious: a proof is a series of logical steps based on some axioms and deduction rules which reaches a desired conclusion. Every step in a proof can be checked for correctness by examining it to ensure that it is logically sound. In David Hilbert's words: "The rules should be so clear, that if somebody gives you what they claim is a proof, there is a mechanical procedure that will check whether the proof is correct or not, whether it obeys the rules or not."

They note, however, that:

> In 1976, Kenneth Appel and Wolfgang Haken proved the 4CT (*Four Colour Theorem*) ... No human being could ever actually read the entire proof to check its

> correctness ... "The real question is this: If no human being can ever hope to check a proof, is it really a proof ?"

The authors' perception of the relation between truth and provability is reflected in their comments:

> ... Godel's incompleteness theorem (GIT) which says that every formal system which is (1) finitely specified, (2) rich enough to include the arithmetic, and (3) consistent, is incomplete[1]. That is, there exists an arithmetical statement which (A) can be expressed in the formal system, (B) is true, but (C) is unprovable within the formal system. ... But what does it mean to be a "true arithmetical statement"? It is a statement about non-negative integers which cannot be invalidated by finding any combination of non-negative integers that contradicts it. ... a true arithmetical statement is a "primordial mathematical reality". ... The essence of GIT is to distinguish between truth and provability. A closer analogy in real life is the distinction between truths and judicial decisions, between what is true and what can be proved in court. How large is the set of true and unprovable statements? If we fix a formal system satisfying all three conditions in GIT, then the set of true and unprovable statements is topologically "large" (constructively, a set of second Baire category, and in some cases even "larger").

Prima facie, under the standard interpretations of classical mathematical theory, which the authors [Ca01] seem to implicitly assume, the above remarks can be taken to imply that the authors accept mathematical truth as being unverifiable effectively. It follows

---

[1] We note that this is an instance of the standard interpretation of Gödel's seminal 1931 paper [Go31a] that may be misleading; we show in Anand [An02b] that Gödel's meta-proof of his Theorem VI in the cited paper actually establishes his formal system P as omega-inconsistent. Hence, Theorem VI holds vacuously, and P is not incomplete.

that there could, then, be any number of (equally reasonable) ways of responding to their question:

> ... what does it mean to be a "true arithmetical statement"?

However, in their attempt to offer an ambitious interpretation of classical theory, the authors do not address the question:

> Can such latitude in the perception of fundamental meta-mathematical concepts such as truth and provability reflect a basic ambiguity in our definitions of foundational mathematical concepts?

On the contrary, the authors seem to be comfortable with the, implicitly Platonic, suggestion that classical concepts of mathematical proof, and even truth, might actually lie beyond the ambit of direct intuitive cognition! They conclude that [Ca01]:

> If we accept the above assumptions about the biological and physical nature of proofs, then there is little 'intrinsic' difference between traditional and 'unconventional' types of proofs as i) first and foremost, we have not access to truth, ii) correctness is not absolute, but nearly certain as mathematics advances by making mistakes and correcting and re–correcting them ..., iii) non–deterministic and probabilistic proofs do not allow mistakes in the applications of rules, they are just indirect forms of checking ... which correspond to various degrees of rigour, iv) the explanatory component, the understanding 'generated' by proofs, while extremely important from a cognitive point of view, is subjective and has no bearing on formal correctness.

> ... more research will be performed in large computational environments where we might or might not be able to determine what the system has done or why ... The blend of logical and empirical-experimental arguments are here to stay and develop.

... There are many reasons which support this prediction. They range from economical ones (powerful computers will be more and more accessible to more and more people), social ones (sceptical oldsters are replaced naturally by youngsters born with the new technology, results and success inspire emulation) to pure mathematical (new challenging problems, wider perspective) and philosophical ones (note that incompleteness is based on the analysis of the computer's behaviour).

## 2. Interpreting classical mathematical theory

### 2.1 Standard interpretations of foundational concepts may be ambiguous

Now, in [An03e] we note that, as is implicit in Mendelson's [Me90] following remarks (*italicised parenthetical qualifications added*), standard interpretations of classical foundational concepts can, indeed, be argued as being either ambiguous, or non-constructive, or both:

> Here is the main conclusion I wish to draw: it is completely unwarranted to say that CT (*Church's Thesis*) is unprovable just because it states an equivalence between a vague, imprecise notion (effectively computable function) and a precise mathematical notion (partial-recursive function). ... The concepts and assumptions that support the notion of partial-recursive function are, in an essential way, no less vague and imprecise (*non-constructive, and intuitionistically objectionable*) than the notion of effectively computable function; the former are just more familiar and are part of a respectable theory with connections to other parts of logic and mathematics. (The notion of effectively computable function could have been incorporated into an axiomatic presentation of classical mathematics, but the acceptance of CT made this unnecessary.) ... Functions are defined in terms of sets, but the concept of set is no clearer (*not more non-constructive, and intuitionistically objectionable*) than that of function and a foundation of mathematics can be based on a theory using function as primitive notion instead of set. Tarski's definition of truth is formulated in set-theoretic terms, but the notion of set is no clearer (*not more non-constructive, and intuitionistically objectionable*), than that of truth. The model-theoretic definition of logical validity is based ultimately on set theory, the foundations of which are no clearer (*not more non-constructive, and intuitionistically objectionable*) than our intuitive (*non-constructive, and intuitionistically objectionable*) understanding of logical validity. ... The notion of Turing-computable function is no clearer (*not more non-constructive, and intuitionistically objectionable*) than, nor more mathematically useful (foundationally speaking) than, the notion of an effectively computable function.

The questions thus arise: Could the thesis conjectured in [Ca01] also be founded on ambiguities that are rooted in the standard interpretations of classical foundational concepts such as "mathematical object", "effective computability", "truth of a formula under an interpretation", "set", "Church's Thesis" etc.; ambiguities that may, moreover, encourage non-constructive, Platonic, interpretations by default? How would such a thesis fare if we could make these concepts unambiguous, and constructive, in an intuitionistically unobjectionable way?

**2.2  Can classical concepts be defined constructively?**

In Anand [An02c], we argue that we can, indeed, define these concepts constructively in terms of a small number of primitive, formally undefined but intuitively unobjectionable, mathematical terms as below:

8(*i*) **Primitive mathematical object**: A primitive mathematical object is any symbol for an individual constant, predicate letter, or a function letter, which is defined as a primitive symbol of a formal mathematical language.[2]

(*ii*) **Formal mathematical object**: A formal mathematical object is any symbol for an individual constant, predicate letter, or a function letter that is either a primitive mathematical object, or that can be introduced through definition into a formal mathematical language without inviting inconsistency.[3]

(*iii*) **Mathematical object**: A mathematical object is any symbol that is either a primitive mathematical object, or a formal mathematical object.

(*iv*) **Set**: A set is the range of any function whose function letter is a mathematical object.

(*v*) **Individual computability**: A number-theoretic function $F(x)$ is individually computable if, and only if, given any natural number $k$, there is an individually effective method (which may depend on the value $k$) to compute $F(k)$.

(*vi*) **Uniform computability**: A number-theoretic function $F(x)$ is uniformly computable if, and only if, there is a uniformly effective method (necessarily independent of $x$) such that, given any natural number $k$, it can compute $F(k)$.

(*vii*) **Effective computability**: A number-theoretic function is effectively computable if, and only if, it is either individually computable, or it is uniformly computable.[4]

---

[2] We note that, as remarked by Mendelson [Me90], the terms "function" and "function letter" - and, presumably, "individual constant", "predicate", and "predicate letter" - can be taken as undefined, primitive foundational concepts.

[3] We highlight the significance of this definition in Meta-lemma 1 in Anand [An02c].

[4] We note that classical definitions of the effective computability of a function (cf. [Me64], p207) do not distinguish between the two cases. The standard interpretation of effective computability is to implicitly

(*viii*) **Individual truth**: A string $[F(x)]$[5] of a formal system P is individually true under an interpretation M of P if, and only if, given any value $k$ in M, there is an individually effective method (which may depend on the value $k$) to determine that the interpreted proposition $F(k)$ is satisfied in M.[6]

(*ix*) **Uniform truth**: A string $[F(x)]$ of a formal system P is uniformly true under an interpretation M of P if, and only if, there is a uniformly effective method (necessarily independent of $x$) such that, given any value $k$ in M, it can determine that the interpreted proposition $F(k)$ is satisfied in M.

(*x*) **Effective truth**: A string $[F(x)]$ of a formal system P is effectively true under an interpretation M of P if, and only if, it is either individually true in M, or it is uniformly true in M.[7]

(*xi*) **Individual Church Thesis**: If, for a given relation $R(x)$, and any element $k$ in some interpretation M of a formal system P, there is an individually effective method such that it will determine whether $R(k)$ holds in M or not, then every element of the domain D of M is the interpretation of some term of P, and there is some P-formula $[R'(x)]$ such that:

---

treat it as equivalent to the assertion: A number-theoretic function $F(x)$] of a formal system P is effectively computable if, and only if, it is both individually computable, and uniformly computable.

[5] We use square brackets to distinguish between the uninterpreted string $[F]$ of a formal system, and the symbolic expression "$F$" that corresponds to it under a given interpretation that unambiguously assigns formal, or intuitive, meanings to each individual symbol of the expression "$F$".

[6] In Anand [An02c], we argue that, under a constructive interpretation of formal Peano Arithmetic, Gödel's undecidable proposition, is individually, but not uniformly, true under the standard interpretation. See also Anand [An03d].

[7] We note that, classically, Tarski's definition of the truth of a formal proposition under an interpretation (cf. [Me64], p49-52) does not distinguish between the two cases. The implicitly accepted (standard) interpretation of the definition appears, prima facie, to be the non-constructive assertion: A string $[F(x)]$ of a formal system P is true under an interpretation M of P if, and only if, it is both uniformly true in M, and individually true in M.

$R(k)$ holds in M if, and only if, $[R'(k)]$ is P-provable.

(In other words, the Individual Church Thesis postulates that, if a relation $R$ is effectively decidable individually (possibly non-algorithmically) in an interpretation M of some formal system P, then $R$ is expressible in P, and its domain necessarily consists of only mathematical objects, even if the predicate letter $R$ is not, itself, a mathematical object.)

(*xii*) **Uniform Church Thesis**: If, in some interpretation M of a formal system P, there is a uniformly effective method such that, for a given relation $R(x)$, and any element $k$ in M, it will determine whether $R(k)$ holds in M or not, then $R(x)$ is the interpretation in M of a P-formula $[R(x)]$, and:

$R(k)$ holds in M if, and only if, $[R(k)]$ is P-provable.

(Thus, the Uniform Church Thesis postulates that, if a relation $R$ is effectively decidable uniformly (necessarily algorithmically) in an interpretation M of a formal system P, then, firstly, $R$ is expressible in P, and, secondly, the predicate letter $R$, and all the elements in the domain of the relation $R$, are necessarily mathematical objects.)

## 2.3  Standard interpretations may admit ambiguity

As we remark in Anand [An03f], the vagueness, and implicitly implied non-constructivity, alluded to in Mendelson's remarks may simply:

> ... reflect, and result from, the non-specification of an effective method for determining that the infinity of intuitive assertions, which are implicit in Tarski's definition of the truth of a formula of a formal system under an interpretation, are, indeed, individually verifiable.

Thus, Tarski's definitions may be seen as implicitly implying, firstly, that relationships may exist only Platonically between the abstract elements of the domain of some model M, since there is no assurance that every such element is constructively definable, or representable, in every model M of PA[8]; and, secondly, that even when such relationships, in some cases, are asserted as holding in M intuitively, this may not be in any effectively verifiable manner.

As we note subsequently [An03f]:

Clearly, what we see here is the thin end of the wedge that keeps the door ajar for the entry of non-constructive, Platonic, elements into the standard interpretations of classical theories; elements that can then be interpreted ambiguously - often creating a mathematical tower of Babel containing frustrated purists, and confused neophytes!

**2.4 Reducing Tarskian truth and satisfiability to provability**

We also argue further that [An03f]:

... if we introduce the concept of effective truth, based on effective methods of verification as suggested above, then we effectively reduce any verifiable truth in the model M to provability in PA. In other words, what such constructive definitions and theses essentially suggest is that, in order to make Tarski's definitions of truth and satisfiability effectively verifiable in any model M of PA, we should be able to argue that:

(*i*) If a string [$R(x)$] is PA-provable, then its interpretation $R(x)$ in M is *both* individually true *and* uniformly true; hence, viewed as a Boolean function, $R(x)$ is Turing-computable.

---

[8] Standard first order Peano Arithmetic such as Mendelson's formal system S ([Me64], p102).

(*ii*) If a string [*R*(*n*)] is PA-provable for any given numeral [*n*], then the interpretation *R*(*x*) in M is *either* individually true *or* uniformly true; hence, viewed as a Boolean function, *R*(*x*) is not necessarily Turing-computable.

(*iii*) If *R*(*x*) is individually true in M (which we may express as (I*x*)*R*(*x*)), then, *R*(*x*) is expressible (cf. [Me64], p117, §2) in PA; hence every element of the domain of M is the interpretation of some term of PA.

(*iv*) *Uniform Turing Thesis*[9]: If *R*(*x*) is uniformly true in M (which we may express as (U*x*)*R*(*x*)), then [*R*(*x*)] is PA-provable; hence, viewed as a Boolean function, *R*(*x*) is Turing-computable.

Now, whilst (*i*), and (*ii*), seem, prima facie, consistent with standard interpretations of Tarski's definitions, (*iii*) clearly does not follow from them; however, as the Löwenheim-Skolem theorem ([Me64], p69, Corollary 2.16) suggests, it may not be inconsistent with such interpretations. It is not obvious whether (*iv*) is independent of, equivalent to, or a consequence of the Individual and Uniform Church Theses.

## 2.5 Some consequences of a constructive interpretation

We then conclude that [An03f]:

> ... the significance of constructively interpreting foundational concepts and assertions of classical mathematics is that:
>
> (*i*) The Uniform Church Thesis implies that a formula [*R*] is P-provable if, and only if, [*R*] is uniformly true in some interpretation M of P.

---

[9] We introduce an equivalent statement of this as an independent Quantum Halting Hypothesis in Anand [An02d], where we show that such a thesis allows us to model a deterministic universe that is essentially unpredictable.

(*ii*) The Uniform Church Thesis implies that if a number-theoretic relation $R(x)$ is uniformly satisfied in some interpretation M of P, then the predicate letter "*R*" is a formal mathematical object in P (i.e. it can be introduced through definition into P without inviting inconsistency).

(*iii*) The Uniform Church Thesis implies that, if a P-formula [*R*] is uniformly true in some interpretation M of P, then [*R*] is uniformly true in every model of P.

(*iv*) The Uniform Church Thesis implies that if a formula [*R*] is not P-provable, but [*R*] is classically true under the standard interpretation, then [*R*] is individually true, but not uniformly true, in the standard model of P.

(*v*) The Uniform Church Thesis implies that Gödel's undecidable sentence GUS is individually true, but not uniformly true, in the standard model of P.[10]

By defining effective computability, both individually and uniformly, along similar lines, we can give a constructive definition of uncomputable number-theoretic functions:

(*vi*) A number-theoretic function $F(x_1, ..., x_n)$ in the standard interpretation M of P is uncomputable if, and only if, it is effectively computable individually, but not effectively computable uniformly.

This, last, removes the mysticism behind the fact that we can constructively define a number-theoretic Halting function that is, paradoxically, Turing-uncomputable.

---

[10] An intriguing consequence of this argument is considered in Appendix 1 of Anand [An03c].

(*vii*) If we assume a Uniform Church Thesis, then every partial recursive number-theoretic function $F(x_1, ..., x_n)$ has a unique constructive extension as a total function.

(*viii*) If we assume a Uniform Church Thesis, then the classical Halting problem is effectively solvable.[11]

(*ix*) If we assume a Uniform Church Thesis, then not every effectively computable function is classically Turing computable (so Turing's Thesis does not, then, hold).

(*x*) If we assume a Uniform Church Thesis, then not every (partially) recursive function is classically Turing-computable.[12]

(*xi*) If we assume a Uniform Church Thesis, then the class P of polynomial-time languages in the P versus NP problem may not define a formal mathematical object.

Further, since a number-theoretic relation is expressible in P if, and only if, it is recursive ([Me64], p142, Corollary 3.29), it follows, as Mendelson argues, that the classical Church's Thesis can, indeed, be viewed as:

(*xii*) *Church's Theorem*: The Individual Church Thesis implies that a number-theoretic function is effectively computable if, and only if, it is recursive[13].

---

[11] This was not included in [An03f].

[12] The classical proof that every (partially) recursive function is classically Turing-computable uses induction over (partial) recursive functions, thus assuming that every such function is a mathematical object; by Meta-lemma 1, such an assumption is invalid.

[13] We note that the classical Church Thesis is the assertion: "A number-theoretic function is effectively computable if, and only if, it is recursive" (cf. [Me64], p227).

**2.6  Defining formal, constructive and Platonic concepts**

Now, what we have essentially argued above is, as also noted in Anand [An03f], that:

> ... a constructive interpretation of classical mathematical concepts recognises that, since Anand ([An02c], Meta-theorem 1) proves constructively that not every, effectively well-defined, classical, mathematical concept is formalisable in its intuitive entirety, we do not need to treat every non-formalisable mathematical concept as necessarily Platonic, and so outside the reach of effective methods. Ipso facto, we can distinguish between:
>
> (*i*) formalisable mathematical concepts (such as arithmetical functions and relations) that are within the ambit of formal effective methods (provability),
>
> (*ii*) non-formalisable mathematical concepts (such as recursive number-theoretic functions and relations) that are within the ambit of non-formal effective methods (effective truth), and
>
> (*iii*) mathematical concepts that are essentially unverifiable by any effective method, and which implicitly assume the existence of a non-constructive Platonic oracle of intuitive truth.

## 3. Mathematical proof and non-algorithmic effective methods

Accordingly, a more appropriate definition of mathematics and mathematical proof - which can be seen, prima facie, as a unifying thread in all the eight stages sought to be identified as distinct by the authors in [Ca01] - would be:

> Mathematics is a language where proof is the yardstick for unambiguous expression and communication.

## 3.1 Non-algorithmic effective methods: Gödel oracles

The raison d'etre, and possible significance, of such a paradigm shift - from seeing mathematics as an expression of relations between a universe of abstract objects, to viewing it solely as a language of unambiguous communication - is reflected in Anand [An03c] and, more particularly, in the following remarks, reproduced from the author's response to a correspondent[14]:

> Broadly speaking, ... Gödel's meta-reasoning (which is constructive and intuitionistically unobjectionable), actually establishes the existence of an undecidable sentence by means of an effective, and not simply theoretical, Turing-oracle.
>
> If we ignore the details of his proof, what Gödel's oracle does is to effectively prove, meta-mathematically, that Gödel's undecidable arithmetical relation $R(x)$[15] is such that, for any given natural number $n$:
>
> (*i*) PA proves: $[R(n)]$,[16]
>
> but:
>
> (*ii*) PA does not prove: $[(Ax)R(x)]$.[17]

---

[14] Reproduced from the correspondence between Michael N. Christoff and the author in the Yahoo! Group comp-sci-theory dated 1st May 2003.

[15] For "relation $R(x)$", read "proposition $(Ax)R(x)$".

[16] This is actually an implicit meta-lemma in Theorem VI of Gödel's seminal 1931 paper [Go31a]; for an explicit proof, see Anand [An02a], §4.3(*b*)(*viii*).

[17] This, too, is an implicit meta-lemma in Theorem VI of Gödel's seminal 1931 paper [Go31a]; for an explicit proof, see Anand [An02a], §4.3(*a*)(*ix*).

Thus Gödel's oracle effectively establishes, by meta-theorem[18] (*i*), that, given any natural number *n*, we can always find some proof (which can be converted into an effective method) that the string [*R*(*n*)] is PA-provable (even though it gives no clue as to how we should go about finding such a proof in any, individual, case). Hence, the arithmetical sentence *R*(*n*) can be effectively asserted as true, for any natural number *n*, by Tarski's definition of the "truth" of the arithmetical predicate *R*(*x*) under the standard interpretation M of PA.

However, because of meta-theorem (*ii*), there is no effective method for determining a PA-proof of the string [(A*x*)*R*(*x*)], since such a proof does not exist. Hence, there is no guarantee of a uniform effective method (algorithm / Turing machine) such that, given any natural number *n*, it will determine that the sentence *R*(*n*) is true under the standard interpretation M. Such a uniform effective method would, of course, be guaranteed if [(A*x*)*R*(*x*)] were PA-provable.

Now, my intended thesis is that meta-theorem (*ii*), in fact, implies that there is no uniform effective method (algorithm / Turing machine) such that, given any natural number *n*, it will determine that the sentence *R*(*n*) is true under the standard interpretation.

In other words, I intend replacing the *classical Turing Thesis*:

> (*iii*) A number-theoretic predicate *F*(*x*), viewed as a Boolean function, is effectively computable if, and only if, it is classically Turing-computable,

with, for instance, a *Uniform Turing Thesis*:

---

[18] For "meta-theorem", read "meta-lemma" in this section.

(*iv*) A PA-string [(A*x*)*F*(*x*)] is provable if, and only if, in it's standard interpretation M, the predicate *F*(*x*), viewed as a Boolean function, is classically Turing-computable as true for every input.

Now, like the classical Turing Thesis, the Uniform Turing Thesis cannot be proven. However, unlike the classical Turing Thesis, this thesis can never be disproved.

The reason: The unprovability of the PA-string [(A*x*)*R*(*x*)] means that we could never formally prove that a Turing machine, (algorithm) T(*R*), that computes the arithmetical predicate *R*(*x*), when seen as a Boolean function, will halt and return true on every input. In fact, in view of meta-theorems (*i*) and (*ii*), the Uniform Turing Thesis implies that any such T(*R*) will loop for some natural number $k$, even though *R*($k$) is true.

The significance of this is that a Uniform Turing Thesis also implies that the standard Turing Thesis does not hold, since Gödel's oracle is, then, an effective (non-algorithmic) method that is not classically Turing-computable. By Occam's dictum, and since there is no loss of generality in replacing the classical Turing Thesis with a Uniform Turing Thesis, the replacement is to be preferred.

The above is, essentially, the argument for introducing constructive Gödel oracles, so that we may extend the scope of effective methods to include non-algorithmic effective methods by means of a Uniform Turing Thesis.

**3.2 Defining Tarskian truth verifiably**

The significance of replacing the classical Church Thesis by Individual and Uniform Church Theses lies in the fact that, without these, there is no explicit convention for asserting Tarskian truth unambiguously. As we argue in the above correspondence:

One reason why Turing oracles may, so far, be viewed as offering theoretical, rather than effective, solutions to constructive problems of Theoretical Computer Science (such as, say, the Halting problem), could be that they are deeply wedded to standard interpretations of classical theory that, in turn, are rooted in Tarski's definitions, such as:

(*i*) The PA-string $[(Ax)F(x)]$ is true under an interpretation M if, and only if, $F(x)$ is satisfied by every element $s$ of M (i.e. $F(s)$ holds in M).

Now, the three significant points to note here are that Tarski implicitly assumes:

(*ii*) there is some language, accessible to us, in which $F(s)$ is expressible for any $s$ of M;

(*iii*) this is a language in which we can verify that $F(s)$ holds for any $s$ of M;

(*iv*) there is some effective method to verify that $F(s)$ holds for any element of M.

Now it follows that, if we seriously intend to make our mathematical language precise and unambiguous, then:

(*v*) we should introduce these assumptions explicitly as premises in Tarski's definitions;

(*vi*) we should specify what we mean by an effective method that can verify that $F(s)$ holds for any element $s$ of M.

Now, if we are not uncomfortable with accepting the Church-Turing Thesis (which is widely assumed to hold), then it is reasonable to assume that whatever we consider as an effective method must necessarily be arithmetic. It then follows that:

> (*vii*) any effective method must either be a uniform effective method (algorithm / Turing machine), or an individual effective method (non-algorithmic / Gödel oracle);
>
> (*viii*) we can only verify that *F*(*s*) holds in M by an effective method if every *s* in M is the interpretation of some symbol [*s*] of PA.

This is, essentially, the reasoning behind my introduction of constructive definitions of classical concepts such as "mathematical object", "individual/uniform effective methods", "individually true/computable predicates/functions", "Individual/Uniform Church Thesis", etc.

Prima facie, it seems to me that the Uniform Turing Thesis, as expressed earlier, is actually a theorem that follows from the introduction of these definitions. However, it may, in fact, be an additional thesis; I am, at the moment, not quite clear whether the Individual / Uniform Church Theses should be treated as definitions or theses.

Note that, if we leave (*ii*), (*iii*), and (*iv*) as implicit, hence undecided and ambiguous, then we are immediately prevented from making constructive assertions that would invalidate non-constructive interpretations of PA. We could, then, develop non-constructive interpretations that are not necessarily claimed as valid or meaningful, but, by default, simply as valid, and possibly meaningful, till shown otherwise.

However ... if we tolerate such ambiguities, then we are in danger of shortchanging scientific disciplines for whom mathematics is, essentially, a language of reliable, and verifiable, external expression and communication. Such a language should, clearly, be based on notions of formal truth that offer a maximum of precision in, and verifiability of, its assertions, with a minimum of ambiguity.

> For most scientific disciplines, the authority of the standard interpretations of classical mathematics is seen, and accepted - perhaps with some element of reluctance, since such acceptance occasionally flies against the grain of observation and experience - not only as absolute, but also as implicitly promising sufficiency, when needed, to help bridge the seemingly unbridgeable chasm between a Platonic world of abstract objects, and the real world of sensory perceptions, that sometimes confronts such disciplines!

## 4. The Halting probability

### 4.1 CCM's Thesis and the Halting problem

The significance of a constructive interpretation of classical concepts, as outlined above, for the arguments offered by [Ca01], emerges if we note that the:

> *CCM Thesis*: Turing-computability is equivalent to PA-provability.

is implicit in the following remarks [Ca01]:

> Classically, there are two equivalent ways to look at the mathematical notion of proof: logical, as a finite sequence of sentences strictly obeying some axioms and inference rules, and computational, as a specific type of computation. Indeed, from a proof given as a sequence of sentences one can easily construct a Turing machine producing that sequence as the result of some finite computation and, conversely, given a machine computing a proof we can just print all sentences produced during the computation and arrange them into a sequence.

Now, prima facie, the *CCM Thesis* seems equivalent to the Uniform Church Thesis, since it can be interpreted as implying that if a function $f(x)$ is uniformly (Turing) computable, then, firstly, *f* is a mathematical object and, secondly, the formula [(E!*y*)*f*(*x*) = *y*] is PA-provable.

It would then follow, from Meta-lemma 14 in Anand [An02c], that the Halting problem is effectively solvable; for, what we essentially argue there is that, given any Turing machine T, if there is an effective method to recognise whether or not a given string P is a valid input of T, then a *CCM Thesis* should imply that there is an individual, non-algorithmic, effective method that will determine whether or not T will halt on input P.

**4.2 An effective solution of the Halting problem**

We reproduce this argument below.

**Theorem 1**: The *CCM Thesis* implies that the Halting problem is effectively solvable.

**Proof**: Given a Turing machine that computes a number-theoretic function *F*(*x*), we note that, by ([Me64], p233, Corollary 5.13), *F*(*x*) is partial recursive. We may thus assume that such an *F* is obtained from a recursive function *G* by means of the unrestricted $\mu$-operator; in other words, that (cf. [Me64], p214):

$F(x) = \mu y(G(x, y) = 0)$.

If [*H*(*x*, *y*)] expresses ~(*G*(*x*, *y*) = 0) in a formal system of Arithmetic such as standard PA, we then consider the PA-provability, and truth in the standard interpretation M of PA, of the formula [*H*(*a*, *y*)] for a given numeral [*a*] of PA, as below:

(*a*) Let Q1 be the meta-assertion that [*H*(*a*, *y*)] is not effectively true in M. Hence there is no effective method in M to determine that, for any given *y* in M, *y* satisfies [*H*(*a*, *y*)] classically. It follows that there is no uniformly effective method (algorithm / Turing machine) in M to determine that, for any given *y* in M, *y* satisfies [*H*(*a*, *y*)] classically.

Since $G(a, y)$ is recursive, it follows that there is some finite $k$ such that any Turing machine T1($y$) that computes $G(a, y)$ will halt and return the value 0 for $y = k$.

(*b*) Next, let Q2 be the meta-assertion that [$H(a, y)$] is effectively true in the standard interpretation M of PA, but that there is no uniformly effective method (algorithm / Turing machine) in M to determine that, for any given $y$ in M, $y$ satisfies [$H(a, y)$] classically.

Since $G(a, y)$ is recursive, it follows that there is some finite $k$ such that the T1($y$) will halt, and return a symbol for self-termination (looping[19]) for $y = k$.

(*c*) Finally, let Q3 be the meta-assertion that [$H(a, y)$] is effectively true in the standard interpretation M of PA, and that there is a uniformly effective method in M to determine that, for any given $y$ in M, $y$ satisfies [$H(a, y)$] classically. We then have that that [$H(a, y)$] is uniformly true in the standard interpretation M of PA.

If we assume a *CCM Thesis*, it then follows (cf. [An02], Meta-lemma 8) that [$H(a, y)$] is PA-provable. Let $h$ be the Gödel-number of [$H(a, y)$]. We consider, now, Gödel's primitive recursive number-theoretic relation $x$B$y$ ([Go31a], p22, definition 45), which holds in M if, and only if, $x$ is the Gödel-number of a proof sequence in PA for the PA-formula whose Gödel-number is $y$. It follows that there is some finite $k$ such that any Turing machine T2($y$), which computes the characteristic function of $x$B$h$, will halt and return the value 0 for $x = k$.

---

[19] We note that any Turing machine can be designed to recognise a looping situation; it simply records every instantaneous tape description at the execution of each machine instruction, and compares the current instantaneous tape description with the record. It can thus be meta-programmed to abort a loop, and return a meta-symbol indicating self-termination.

Since Q1, Q2, and Q3 are mutually exclusive and exhaustive, it follows that, when run simultaneously over the sequence 1, 2, 3, ... of values for $y$, one of {T1(y) // T2(y)} will always halt for some finite value of $y$ for any given $a$.

Thus, the Halting problem is effectively solvable if we assume a *CCM Thesis*.

**4.3 Is the Halting probability a Dedekind real number?**

Now, the possible digital strings of length $k$ (without repetitions) are less than $2^k$. If $f(k)$ of these, when input to a given Turing machine T, yield Halting programs, the classical probability that any given string of length $k$ will halt is greater than $HP_{T, k} = f(k)/2^k$.

The classical probability that any given string of length less than, or equal to, $k$ will halt is thus greater than $HP_{T, i=<k} = (Sum_{i=<k} f(i))/ (Sum_{i=<k} 2^i)$.

Since, by Theorem 1 above, we can determine $f(i)$ for any given $i$, we can also determine $HP_{T, i=<k}$ as a rational number for any given $k$.

The question arises: Does the non-terminating sequence $HP_T$, expressed by < $HP_{T, i=<1}$, $HP_{T, i=<2}$, $HP_{T, i=<3}$, ...> define a Dedekind cut?

Now, since $f(k)$ is determined by a parallel duo of Turing machines that is not, itself, a Turing machine, but which can be viewed as a deterministic Turing oracle, it follows that there is no uniformly effective method of determining $f(k)$ (cf. [An02c], Corollary 14.2); in other words, although the number-theoretic function $f(x)$ is a well-defined mathematical concept, $f$ is not a mathematical object. Hence $f(x)$ may be considered as determinate, but uncomputable; its values are essentially unpredictable, and so, by definition, truly random.

It follows that the rational function $HP_{T,\ i=<k}$ (viewed as a number-theoretic function over ordered integer pairs), too, is, essentially, (Turing) uncomputable. Hence $HP_T$ is also not a mathematical object, and so its range does not define a Dedekind real number.

Moreover, since $HP_{T,\ i=<k}$ is effectively computable individually for any given natural number $k$, we may look upon $HP_T$ as a random mathematical function that, nevertheless, contains an infinitude of non-algorithmically computable information.

**4.4 Standard interpretations of the significance of the Halting probability**

We contrast the above with the standard interpretation of the relation of classical theory to the Halting problem, and of the significance of the Halting probability, as expressed by the authors in [Ca01]:

> In modern times a penetrating insight into the incompleteness[20] phenomenon has been obtained by Chaitin's information–theoretic analysis ... The simplest way to state one of Chaitin's main results is the following ... : "If you have N bits of axioms, then you can never prove that a program is the smallest possible if it is more than N bits long." The most striking results have been obtained by studying the Chaitin's Omega Number, *Omega*, the halting probability of a self-delimiting universal Turing machine. This number is not only uncomputable, but also (algorithmically) random. Chaitin has proven the following important theorem: If ZFC is arithmetically sound[21], then, ZFC can determine the value of only finitely many bits of *Omega*, and one can give a bound on the number of bits of *Omega* which ZFC can determine.[22] Robert

---

[20] The significance of this may need to be viewed, however, in the light of earlier remarks in footnote 1.

[21] This theorem would, then, hold vacuously; it follows from Anand ([An02c], Corollary 1.1) that ZFC is not arithmetically sound, in the sense that it is not a model for standard PA..

[22] This result may also need to be viewed against the arguments of the previous section; if the Halting probability is not a Dedekind real number, then it is not clear what arithmetical interpretation, or significance, is to be given to a routine that calculates "finitely many bits of *Omega*" or one that "can give a bound on the number of bits of *Omega* which ZFC can determine".

> Solovay ... has constructed a self-delimiting universal Turing machine such that ZFC, if arithmetically sound, cannot determine any single bit of its halting probability ... Rephrased, the most powerful formal axiomatic system is powerless when dealing with the questions of the form "is the $m$'th bit of *Omega* 0?" or "is the $m$'th bit of *Omega* 1?".

Clearly, the significant difference is that, since standard interpretations of classical theory do not explicitly distinguish between well-defined mathematical concepts and well-defined mathematical objects, the distinction between a non-constructively defined Chaitin real number (*Omega*), and a Dedekind real number (which provides the foundation for classical mathematical theory) remains implicit; it is thus not clear whether the significance ascribed to mathematical assertions about the former can be meaningfully, and significantly, extended to the latter[23].

## 5. Are mechanistic proofs of mathematical problems logically sound?

In conclusion, we note that the pedantic point made by Robertson, Sanders, Seymour and Thomas, as quoted in [Ca01], about the reliability of their 1996 computer-generated proof of 4CT, can be seen as besieged sophistry:

> However, an argument can be made that our "proof" is not a proof in the traditional sense, because it contains steps that can never be verified by humans. In particular, we have not proved the correctness of the compiler we compiled our programs on, nor have we proved the infallibility of the hardware we ran our programs on. These have to be taken on faith, and are conceivably a source of error.

---

[23] In Anand [An03b], we argue that standard interpretations of classical mathematical theory may be held responsible for (misleadingly) encouraging a similar blurring of the distinction between Cantorian real numbers and Dedekind real numbers.

Although one may doubt whether their program is logically sound, there is no essential reason why such soundness cannot be established theoretically[24]; thus, there is no sensible reason to doubt the output, even though the relation between a logically effective method and a mechanistic computation is, indeed, one of faith.

It is conceivable that an appropriately designed, and maintained, machine continuously calculating the digits of Pi may start outputting an unending series of zero's, perhaps for as long as a zillion years. Our belief that it will eventually output a digit other than zero comes about because of our faith that the machine is faithfully translating our theoretical calculations concerning the digits of Pi into a physical language of mechanical I/O devices; hence we believe, as an article of faith, and despite the physical evidence, that the series of zeros will not be unending.

Similarly, the effective solution of the Halting problem consists of a duo of parallel Turing machines, where an assumption of a Uniform Church (or Turing / CCM) Thesis implies that one of the two will halt either because its program is a halting program, or because it recognises and halts due to an impending looping situation. Again, physically, this could take a zillion years in some particular case, and may call upon all our reserves of patience and faith.

So, finally, what we should place our faith in is our ability to intuitively "see" the soundness of our axiomatic assertions, and theses, at the lowest, and not, as some believe, at a sufficiently sophisticated, level of understanding. This, of course, is the distinguishing feature of, for instance, Dedekind's formulation of the Peano axioms, or - as any high-school student can testify - Euclid's axioms for geometry; except the parallel postulate, these axioms are obvious to all, even if their consequences are not.

---

[24] Since every step of a formal proof sequence is either an axiom, or an immediate consequence of any two preceding elements of the sequence, each step can be effectively verified mechanistically by identifying the concerned axiom, or two preceding elements of the sequence. So long as each step is verified as logically sound, such a procedure need not be time bound, nor limited to the conceptual ability of any one individual to grasp, or even verify, the correctness of the entire proof.